\documentclass[a4paper]{article}

\usepackage{amsmath,amssymb,amsthm}
\usepackage{graphicx}
\usepackage{subcaption}
\graphicspath{{./figures/}}

\theoremstyle{definition}
\newtheorem{theorem}{Theorem}
\newtheorem{definition}[theorem]{Definition}
\newtheorem{lemma}[theorem]{Lemma}
\newtheorem{remark}[theorem]{Remark}
\newtheorem{corollary}[theorem]{Corollary}

\usepackage{xifthen}
\newcommand{\textcite}[2][]{\ifthenelse{\isempty{#1}}{\cite{#2}}{\cite[#1]{#2}}}
\newcommand{\parencite}[2][]{\ifthenelse{\isempty{#1}}{\cite{#2}}{\cite[#1]{#2}}}

\begin{document}

\title{Statistical tests based on R\'{e}nyi entropy estimation}
\author{Mehmet Siddik Cadirci$^1$, Dafydd Evans$^1$, Nikolai Leonenko$^1$ \\ and Oleg Seleznjev$^2$}
\date{
	\slshape\small 
	$^1$ School of Mathematics, Cardiff University, Cardiff, Wales, UK. \\
	$^2$ Department of Mathematics and Mathematical Statistics, Ume\aa\ University, \textsl{Ume\aa }, Sweden.\\[3ex] 
	\normalfont
	\today 
}
\maketitle

\begin{abstract}
Entropy and its various generalizations are important in many fields, including mathematical statistics, communication theory, physics and computer science, for characterizing the amount of information associated with a probability distribution. In this paper we propose goodness-of-fit statistics for the multivariate Student and multivariate Pearson type II distributions, based on the maximum entropy principle and a class of estimators for R\'{e}nyi entropy based on nearest neighbour distances. We prove the $L^2$-consistency of these statistics using results on the subadditivity of Euclidean functionals on nearest neighbour graphs, and investigate their rate of convergence and asymptotic distribution using Monte Carlo methods. In addition we present a novel iterative method for estimating the shape parameter of the multivariate Student and multivariate Pearson type II distributions.


\end{abstract}

\section{Introduction}
\label{sec:intro}

Entropy is a measure of randomness that emerged from information theory, and its estimation plays an important role in many fields including mathematical statistics, cryptography, machine learning and indeed almost every branch of science and engineering. There are many possible definitions of entropy, for example, the \emph{differential entropy} of a multivariate density function $f:\mathbb{R}^m\to\mathbb{R}$ is defined by
\begin{equation}\label{eq:differential_entropy}
H_1(f) = -\int_{\mathbb{R}^{m}}f(x)\log f(x)\,dx.
\end{equation}

In this paper, we propose statistical tests for a class of multivariate Student and Pearson type II distribtuions, based on estimation of their R\'{e}nyi entropy
\begin{equation}\label{eq:renyi_entropy}
H_q(f) = \frac{1}{1-q}\log\int_{\mathbb{R}^m}f^q(x)\,dx,\quad q\neq 1.
\end{equation}

Estimation of R\'{e}nyi entropy for absolutely continuous multivariate distributions has been considered by many authors, including 
\textcite{kozachenko1987},
\textcite{goria2005},
\textcite{evans2008},
\textcite{leonenko2008},
\textcite{leonenko2010-correction},
\textcite{penrose2011}, 
\textcite{delattre2017},
\textcite{gao2018}, 
\textcite{bulinski2019}, and
\textcite{berrett2019}.
The quadratic R\'{e}nyi entropy was investigated by \textcite{leonenko2010a}.
An entropy-based goodness-of-fit test for generalized Gaussian distributions is presented by \textcite{cadirci2020}.
A recent application to image processing can be found in \textcite{dresvyanskiy2020}.

The remainder of this paper is organized as follows. 
In Section~\ref{sec:maxent}, we present maximum entropy principles for R\'enyi entropy.
In Section~\ref{sec:eest}, we provide nearest-neighbour estimators for R\'enyi entropy.
In Section~\ref{sec:htest}, we propose statistical tests for the multivariate Student and Pearson II distributions.
In Section~\ref{sec:numer}, we report the results of numerical experiments.

\section{Maximum entropy principles}
\label{sec:maxent}

Let $X\in\mathbb{R}^m$ be a random vector that has a density function $f(x)$ with respect to Lebesgue measure on $\mathbb{R}^m$, and let $S=\{x\in\mathbb{R}^m:f(x)>0\}$ be the support of the distribution.
The R\'{e}nyi entropy of order $q\in(0,1)\cup(1,\infty)$ of the distribution is
\begin{equation}\label{eq:renyi}
H_q(f) = \frac{1}{1-q}\log\int_S f^{q}(x)\,dx,
\end{equation}
which is continuous and non-increasing in $q$. If the support has finite Lebesgue measure $|S|$,  then 
\[
\lim_{q\to 0}H_q(f) = \log|S|,
\]
otherwise $H_q(f)\to\infty$ as $q\to 0$. Note also that
\[
\lim_{q\to 1}H_q(f) = H_1(f) = -\int_{S}f(x)\log f(x)\,dx.
\]

%
Let $a\in\mathbb{R}^m$ and let $\Sigma $ be a symmetric positive definite $m\times m$ matrix. 
%
%
The \emph{multivariate Gaussian} distribution $N_m(a,\Sigma)$ on $\mathbb{R}^m$ has density function
\[
f_{a, \Sigma}^{G}(x) = (2\pi)^{-m/2} |\Sigma|^{-1/2}
\exp\left(-\frac{1}{2}(x-a)^\prime\, \Sigma^{-1}(x-a)\right).
\]
For $X\sim N_m(a,\Sigma)$, we have $a=\mathbb{E}(X)$ and $\Sigma = {\rm Cov}(X)$, where $\text{Cov}(X)=\mathbb{E}[(X-a)(X-a)']$ is the covariance matrix of the distribution.

\smallskip
For $\nu>0$, the \emph{multivariate Student distribution} $T_m(a,\Sigma,\nu)$ on $\mathbb{R}^m$ has density function
\begin{align}
& f^S_{a,\Sigma,\nu}(x) = c_S |\Sigma|^{-1/2}
\left(1 + \frac{1}{\nu}(x - a)'\Sigma^{-1}(x-a)\right)^{-\frac{\nu+m}{2}} \nonumber\\
&\text{where } c_S(m,\nu) = \frac{\Gamma[(\nu+m)/2]}{(\pi\nu)^{m/2}\Gamma(\nu/2)}.\label{eq:const-student}
\end{align}
For $X\sim T_m(a,\Sigma,\nu)$ we have $a=\mathbb{E}(X)$ when $\nu>1$ and $\Sigma = (1-2/\nu)\text{Cov}(X)$ when $\nu>2$, see \parencite{johnson2007}.
It is known that $f^S_{a,\Sigma,\nu}(x)\to f^G_{a,\Sigma}(x)$ as $\nu\to\infty$.

\smallskip
For $\gamma>0$, the \emph{multivariate Pearson Type II} distribution $P_m(a,\Sigma,\gamma)$ on $\mathbb{R}^m$, also known as the \emph{Barenblatt} distribution, has density function
\begin{align}
& f^P_{a,\Sigma,\gamma}(x) = c_P|\Sigma|^{-1/2}
\left[1-(x-a)^{\prime}\Sigma^{-1}(x-a)\right]^{\gamma}_{+} \nonumber \\
&\text{where } t_+=\max\{t,0\} \text{ and }
c_P(m,\gamma) = \frac{\Gamma(m/2+\gamma+1)}{\pi^{m/2}\Gamma(\gamma+1)}.\label{eq:const-pearsonII}
\end{align}
For $X\sim P_m(a,\Sigma,\gamma)$ we have $a = \mathbb{E}(X)$ and $\Sigma = (m+2\gamma+2)\text{Cov}(X)$.
It is known that $f^P_{a,\Sigma,\gamma}(x)\to f^G_{a,\Sigma}(x)$ as $\gamma\to\infty$.

\begin{remark}
If the covariance matrix $C$ is diagonal, the Pearson Type II distribution belongs to the class of time-dependent distributions 
\[
u(x,t) = c(\beta ,\gamma)t^{-\alpha m}\left(1-\left(\frac{\|x\|}{ct^\alpha}\right)^\beta\right)_{+}^{\gamma}
\]
with $c>0$, $\mathrm{supp}\{u(x,t)\}=\{x\in\mathbb{R}^m:\|x\|\leq ct^\alpha\}$ and
\[
c(\beta ,\gamma) = \beta \Gamma \left( \frac{m}{2}\right) /\left[ 2c^{m}\pi
^{\frac{m}{2}}B\left( \frac{m}{\beta },\gamma +1\right) \right] ,
\]
which are known as \emph{Barenblatt solutions} of the source-type non-linear diffusion equations $u'_t=\Delta(u^q)$, where $q>1$, $\Delta$
is the Laplacian and $\gamma=1/(q-1)$. For details, see \textcite{frank2005}, \textcite{vazquez2007} and \textcite{degregorio2020}.
\end{remark}

\subsection{R\'enyi entropy} 
\label{subsec:renyientropy}

The R\'enyi entropy of the multivariate Gaussian distribution $N_m(a,\Sigma)$ is
\begin{align*}
H_q(f^G_{a,\Sigma})
	& = \log\left[(2\pi)^{m/2}|\Sigma|^{1/2}\right]-\frac{m}{2(1-q)}\log q \\
	& = H_1(f^G_{a,\Sigma})-\frac{m}{2}\left( 1+\frac{\log {q}}{1-q}\right) 
\end{align*}
where 
$H_1(f^G_{a,\Sigma})=\log\left[(2\pi e)^{m/2}|\Sigma|^{1/2}\right]$ is the differential entropy of $N_m(a,\Sigma )$.
From \textcite{zografos2005}, the  R\'enyi entropy of the multivariate Student distribution $T_m(a,\Sigma,\nu)$ is
\begin{equation}\label{eq:renyi-student}
H_q\big(f^S_{a,\Sigma,\nu}\big) = \frac{1}{2}\log|\Sigma| + c'_S(m,q,\nu)
\end{equation}
where
\begin{eqnarray*}
\lefteqn{c'_S(m,q,\nu) =
\frac{1}{1-q}\log\left(\frac{B\left(q\left(\frac{\nu+m}{2}\right)-\frac{m}{2},\frac{m}{2}\right)}{B\left(\frac{\nu}{2},\frac{m}{2}\right)^q}\right) } \\
	&& \qquad\qquad\qquad\qquad\qquad + \frac{m}{2}\log(\pi\nu) - \log\Gamma\left(\frac{m}{2}\right). 
\end{eqnarray*}
Likewise the R\'{e}nyi entropy of the multivariate Pearson Type II distribution $P_m(a,\Sigma,\gamma)$ is
\begin{equation}\label{eq:renyi-pearsonII}
H_q\big(f^P_{a,\Sigma,\gamma}\big) = \frac{1}{2}\log|\Sigma| + c'_P(m,q,\gamma),
\end{equation}
where
\begin{eqnarray*}
\lefteqn{c'_P(m,q,\gamma) =
\frac{1}{1-q}\log\left(\frac{B\left(q\gamma+1, \frac{m}{2}\right)}{B\left(\gamma+1,\frac{m}{2}\right)^q}\right) } \\
	&& \qquad\qquad\qquad\qquad\qquad + \frac{m}{2}\log(\pi) - \log\Gamma\left(\frac{m}{2}\right). 
\end{eqnarray*}

\subsection{Maximum entropy principle}
\label{subsec:maxent}

\begin{definition}\label{def:classK}
Let $\mathcal{K}$ be the class of density functions supported on $\mathbb{R}^m$ and subject to the constraints
\[
\int_{\mathbb{R}^m}xf(x)\,dx = a 
\quad\text{and}\quad
\int_{\mathbb{R}^m}(x-a)(x-a)^{\prime}f(x)\,dx = C
\]
where $a\in\mathbb{R}^m$ and $C$ is a symmetric and positive definite $m\times m$ matrix. 
\end{definition}
It is well-known that the differential entropy $H_1$ is uniquely maximized by the multivariate normal distribution $N_m(a,\Sigma)$, that is
\[
H_1(f) \leq H_1(f^G_{a,\Sigma}) = \log\left[(2\pi e)^{m/2}|\Sigma|^{1/2}\right]
\]
with equality if and only if $f = f^G_{a,\Sigma}$ almost everywhere.
The following result is discussed by
\textcite{kotz2004}, 
\textcite{lutwak2004},
\textcite{heyde2005}, and 
\textcite{johnson2007} among others.

\begin{theorem}[Maximum R\'enyi entropy]
\label{thm:maxent}
\mbox{}\par
\noindent(1) For $m/(m+2) < q < 1$, $H_q(f)$ is uniquely maximized over $\mathcal{K}$ by the multivariate Student distribution $T_m(a,\Sigma,\nu)$ with $\nu = 1/(1-q)-m$ and $\Sigma=(1-2/\nu)C$.

\smallskip
\noindent(2) For $q>1$, $H_q(f)$ is uniquely maximized over $\mathcal{K}$ by the multivariate Pearson Type II distribution $P_m(a,\Sigma,\gamma)$ with $\gamma = 1/(q-1)$ and $\Sigma=(2\gamma+m+2)C$.
\end{theorem}
Applying \eqref{eq:renyi-student} and \eqref{eq:renyi-pearsonII} yields the following.
\begin{corollary}
\mbox{}\par
\noindent(1) For $m/(m+2)<q<1$ the maximum value of $H_q$ is 
\[
H_q^{\rm max} = \frac{1}{2}\log|\Sigma| + c'_S(m,q,\nu)
\]
with $\nu = 1/(1-q)-m$ and $\Sigma=(1-2/\nu)C$.

\noindent(2) For $q>1$ the maximum value of $H_q$ is 
\[
H_q^{\rm max} = \frac{1}{2}\log|\Sigma| + c'_B(m,q,\gamma)
\]
with $\gamma = 1/(q-1)$ and $\Sigma = \left(2\gamma+m+2\right)C$.
\end{corollary}

\section{Statistical estimation of R\'{e}nyi entropy}
\label{sec:eest}

We state some known results on the statistical estimation of R\'{e}nyi entropy due to 
\textcite{leonenko2008}, and
\textcite{penrose2011}.
Extensions of these results can be found in
\textcite{penrose2003},
\textcite{berrett2019},
\textcite{delattre2017},
\textcite{bulinski2019}, and
\textcite{gao2018}.
Let $X\in\mathbb{R}^m$ be a random vector with density function $f$, and let $G_q(f)$ denote the expected value of $f^{q-1}(X)$,
\[
G_q(f) = \mathbb{E}\big[f^{q-1}(X)\big] = \int_{\mathbb{R}^m}f^{q}(x)\,dx
\]
so that $H_q(f) = \frac{1}{1-q}\log G_q(f)$.

\smallskip
Let $X_1,X_2,\ldots,X_N$ be independent random vectors from the distribution of $X$, and for $k\in\mathbb{N}$ with $k<N$, let $\rho_{i,k,N}$ denote the \emph{$k$-nearest neighbour distance} of $X_i$ among the points $X_1,X_2,\ldots,X_N$, defined to be the $k$th order statistic of the $N-1$ distances $\|X_i-X_j\|$ with $j\neq i$, 
\[
\rho_{i,1,N}\leq\rho_{i,2,N}\leq\cdots\leq\rho_{i,N-1,N}.
\]
We estimate the expectation $G_q(f) = \mathbb{E}(f^{q-1})$ by the sample mean
\[
\hat{G}_{k,N,q} = \frac{1}{N}\sum_{i=1}^{N}\left(\zeta_{i,k,N}\right)^{1-q},
\]
where 
$$\zeta_{i,k,N} = (N-1)C_k V_m \rho_{i,k,N}^m
\quad\text{with}\quad
C_k = \displaystyle\left[\frac{\Gamma(k)}{\Gamma(k+1-q)}\right]^{\frac{1}{1-q}}
$$
and
$V_m = \frac{\pi^{\frac{m}{2}}}{\Gamma\left(\frac{m}{2}+ 1\right)}$
is the volume of the unit ball in $\mathbb{R}^m$.

\begin{definition}
For $r>0$, the \emph{$r$-moment} of a density function $f$ is
\[
M_r(f) = \mathbb{E}(\|X\|^r) = \int_{\mathbb{R}^m} \|x\|^r f(x)\,dx,
\]
and the \emph{critical moment} of $f$ is
\[
r_c(f) = \sup\{r>0:M_r(f)<\infty\}
\]
so that $M_r(f)<\infty $ if and only if $r < r_c(f)$.
\end{definition}
The following result was stated without proof in \textcite{leonenko2010-correction}: here we present the proof.

\begin{theorem}\label{thm:main} 
Let $0<q<1$ and $k\geq1$ be fixed.
\begin{enumerate}
\item  If $G_q(f)<\infty$ and    
\begin{align}
& r_c(f) > \frac{m(1-q)}{q},\label{eq:condition-conv-in-mean} \\
&\text{then }\mathbb{E}\left[\hat{G}_{k,N,q}\right]\to G_q(f) \text{ as $N\to\infty$.}
\label{eq:conv-in-mean}
\end{align}

\item If $G_q(f)<\infty$, $q>\frac{1}{2}$ and   
\begin{align}
& r_c(f) > \frac{2m(1-q)}{2q-1},\label{eq:condition-conv-in-mean-square} \\
&\text{then }\mathbb{E}\left[\hat{G}_{k,N,q}-G_q(f)\right]^2\to 0 \quad\text{as $N\to\infty$.}
\label{eq:conv-in-mean-square}
\end{align}
\end{enumerate}

\end{theorem}

\begin{remark}
If $G_q(f) < \infty$ for $q\in\left(1,\frac{k+1}{2}\right)$ then by \parencite{leonenko2008},
\[
\mathbb{E}\left[\hat{G}_{k,N,q} \right] \to G_q(f)
\text{ and }
\mathbb{E}\left[\hat{G}_{k,N,q}-G_q(f)\right]^2 \to 0
\text{ as $N\to\infty$.}
\]
\end{remark}

\begin{remark}
If $G_q(f)<\infty$ for $q\in(0,1)$ and $f(x)=O(\|x\|^{-\beta})$ as $\|x\|\to\infty$ for some $\beta>m$, then $r_c(f)=\beta-m$ and condition~\eqref{eq:condition-conv-in-mean} is automatically satisfied: see \textcite{penrose2011} for a discussion, and counterexamples showing that conditions~\eqref{eq:condition-conv-in-mean} and \eqref{eq:condition-conv-in-mean-square} cannot be omitted in general.
\end{remark}

{\sl Proof of Theorem~\ref{thm:main}.}
Let us write
\[
\hat{G}_{k,N,q}
	= \frac{1}{N}\sum_{i=1}^{N}\left[(N-1)^{1/m}(C_{k}V_{k})^{1/m}\rho_{i,k,N}\right]^{(1-q)m}.
\]

We show that the method proposed by \textcite{penrose2013} for $k=1$ in fact works for any fixed $k\geq 1$.
By Theorem 2.1 of \parencite{penrose2013}, the uniform integrability condition
\begin{equation}
\label{eq:UIC}
\sup_{N}\mathbb{E}\left[ \left\{ ((N-1)(C_{k}V_{k})\rho_{i,k,N-1}^m\right\}^{(1-q)p}\right] <\infty
\end{equation}
for some $p>1$ (statement 1) or some $p>2$ (statement 2) ensures the $L_p$ convergence of $\hat{G}_{k,N,q}$ to $I_{q}$ as $N\rightarrow \infty$.
Because we only need to obtain a bound on left-hand side of \eqref{eq:UIC}, we can use results on the subadditivity of Euclidean functionals defined on the nearest-neighbors graph \parencite{yukich1998}.
We use the following result (Lemma~3.3) from \textcite{penrose2011}, see also \parencite[p.85]{yukich1998}.
\begin{lemma}
\label{lem:penrose}
Let $0<s<m$. 
If $r_{c}(f)>\frac{ms}{m-s}$, then
\begin{align*}
& \sum_{j=1}^{\infty }2^{js}\left[ P(A_{j})\right] ^{\frac{m-s}{m}}<\infty
\quad\text{where}\quad
P(A_j) =\int_{A_{j}}f(x)\,dx \\
& \text{and } A_j = \mathcal{B}(0,2^{j+1})\setminus \mathcal{B}(0,2^j) \text{ for }j=1,2,\ldots 
\end{align*}
with $\mathcal{B}(0,R)=\{x\in\mathbb{R}^{m}:\|x\|\leq R\}$ and $A_0=\mathcal{B}(0,2)$.
\end{lemma}

We continue the proof of Theorem~\ref{thm:main}. Let $b=(1-q)mp$, and note that we can always choose $p$ to ensure  that $0<1-b/m<1.$
By exchangeability,
\begin{eqnarray}
\lefteqn{\mathbb{E}\left[ (N-1)^{1/m}(C_{k}V_{m})^{1/m}\rho _{i,k,N-1}\right]^b}\qquad\qquad \nonumber\\
	& = & \mathbb{E}\left(\frac{1}{N}\sum_{i=1}^{N}\left[(N-1)^{1/m}(C_{k}V_{m})^{1/m}\rho _{i,k,N-1}\right]^b\right)\label{4.3}\nonumber\\
	& = & \frac{(N-1)^{b/m}}{N}(C_{k}V_{m})^{b/m}\mathbb{E}\left(\sum_{i=1}^{N}\rho_{i,k,N-1}^b\right)  \nonumber \\
	& \leq & (C_{k}V_{m})^{b/m}(N-1)^{b/m-1}\mathbb{E}(\mathcal{L}_{k}^{b}(\mathcal{X}_{N})),  \nonumber
\end{eqnarray}
where $\mathcal{X}_N=\left\{X_1,X_2,\ldots,X_N\right\}$, and for any finite point set $\mathcal{X}\subset\mathbb{R}^m$ and $b>0$ we write 
\[
\mathcal{L}_{k}^{b}(\mathcal{X} )=\sum_{x\in \mathcal{X} }\mathcal{D}_{k}^{b}(x,\mathcal{X} ),
\]
where $\mathcal{D}^b_k(x,\mathcal{X})$ denotes the Euclidean distance from $x$ to its $k$-nearest neighbour in the point set $\mathcal{X}\setminus\{x\}$ when ${\rm card}(\mathcal{X})\geq k$; set $\mathcal{D}^b_k(x,\mathcal{X})=0$ if ${\rm card}(\mathcal{X})\leq k$.
The function $\mathcal{X}\mapsto\mathcal{L}_k^b(\mathcal{X})$ satisfies the subadditivity relation
\begin{equation}\label{eq:subadditivity}
\mathcal{L}_{k}^{b}(\mathcal{X} \cap \mathcal{Y})\leq \mathcal{L}_k^b(\mathcal{X}) + \mathcal{L}_k^b(\mathcal{Y}) + U_k t^b
\end{equation}
for all $t>0$ and finite $\mathcal{X}$ and $\mathcal{Y}$ contained in $[0,t]^m$, where $U_k = 2km^{b/2}$, $b>0$. Indeed, if  $\mathcal{X} $ has more than $k$ elements, the $k$-nearest neighbour distances of points in $\mathcal{X} $ can only become smaller when we add some other set $\mathcal{Y}.$ Hence, \eqref{eq:subadditivity} holds with $U_k=0$ if $\mathcal{X}$ and $\mathcal{Y}$ have more than $k$ elements. If $\mathcal{X}$ has  $k$ elements or fewer, then $\mathcal{L}_k^b(\mathcal{X})$ is zero, but when we add the set $\mathcal{Y}$, we gain at most $k$ new edges from points in $\mathcal{X}$ in the nearest neighbours graph, and each of these is of length most $t\sqrt{m}$ (for more details, see \textcite[pp 101-103]{yukich1998}).

Let $s(N)$ be the largest $j\in N$ such that the set $\mathcal{X}_{N}=\{X_1,X_2,\ldots,X_N\}\cap A_j$ is not empty.
Using ideas from \textcite[p.87]{yukich1998} we have that

\[
\mathcal{X}_N\cap\left(\bigcup_{j=0}^{s(N)}A_j\right) = \bigcup_{j=0}^{s(N)}\,(X_N\cap A_j),
\]
and by the subadditivity property,
\begin{eqnarray*}
\lefteqn{\mathcal{L}_k^b(\mathcal{X}_N) \leq \mathcal{L}_k^b\{X_N\cap A_{s(N)}\}} \\
	&& \qquad\qquad + \mathcal{L}_k^b\left(\mathcal{X}_N\cap\left\{\bigcup_{j=0}^{s(N)-1}A_j\right\}\right) 
	+ U_k 2^{(s(N)+1)b}.
\end{eqnarray*}
Applying subadditivity in the same way to the second term on the right yields
\begin{eqnarray*}
\lefteqn{\mathcal{L}_k^b\left(\mathcal{X}_N\cap\left\{\bigcup_{j=0}^{s(N)-1}A_j\right\}\right)
	\leq \mathcal{L}_k^b(\mathcal{X} _{N}\cap A_{s(N)-1})} \\
	&& \qquad\qquad	+ \mathcal{L}_k^b\left(\mathcal{X}_N\cap\left\{\bigcup_{j=0}^{s(N)-2}A_j\right\}\right)
		+ U_k\left(2^{s(N)}\right)^b.
\end{eqnarray*}
Repeatedly applying subadditivity,  we arrive at 
\begin{align}
\mathcal{L}_{k}^{b}\left( X_{1},\ldots ,X_{N}\right) 
	& \leq \sum_{j=0}^{s(N)}\mathcal{L}_k^b(\mathcal{X}_N\cap A_j) + 2^{b+bs(N)}\frac{U_k}{1-2^{-b}} \nonumber \\
	& \leq \sum_{j=0}^{s(N)}\mathcal{L}_k^b(\mathcal{X}_N\cap A_j) + 2^{bs(N)}M_k \nonumber \\
	& \leq \sum_{j=0}^{s(N)}\mathcal{L}_k^b(\mathcal{X}_N\cap A_j) + M_k\max_{1\leq i\leq N}\|X_i\|^b \label{eq:upper-bound}
\end{align}
for some constant $M_k$ depending on $m$, $k$ and $b$.
From \eqref{eq:subadditivity} and \eqref{eq:upper-bound}, we get
\begin{eqnarray}
\lefteqn{\mathbb{E}\left((N-1)^{1/m}(C_{k}V_{m})^{1/m}\rho _{i,k,N-1}\right)^b} \qquad \nonumber\\
	& \leq & (C_{k}V_{m})^{b/m}(N-1)^{b/m-1}\mathbb{E}\left(\sum_{j=0}^{s(N)}\mathcal{L}_k^{b}(\mathcal{X}_N\cap A_j)\right) \nonumber \\
	& & \quad + \,W_k\,\mathbb{E}\left((N-1)^{b/m-1}\max_{1\leq i\leq N}\|X_i\|^b\right)\label{eq:upper-bound2}
\end{eqnarray}
for some constant $W_k$ depending on $m$, $k$ and $b$.
Using Lemma 3.3 of \textcite{yukich1998} we have
\begin{equation}\label{eq:upper-bound3}
{L}_{k}^{b}(\mathcal{X}) \leq L_0({\rm diam}\mathcal{X})^b({\rm card}\mathcal{X})^{1-b/m}
\end{equation}%
for some constant $L_0>0$. Following \textcite{penrose2011}, by Jensen's inequality and the fact that ${\rm diam}(A_j)=2^j$, we obtain from \eqref{eq:upper-bound2} and \eqref{eq:upper-bound3} that
\begin{eqnarray}
\lefteqn{(N-1)^{b/m-1}\mathbb{E}\left(\sum_{j=0}^{s(N)}{L}_{k}^{b}(\mathcal{X}_N\cap A_j)\right)} \nonumber\\
	&&\qquad\qquad \leq {L}_1\sum_{j=0}^{s(N)}2^{jb}\left[ \mathbb{P}(X_{1}\in A_{j})\right] ^{1-b/m} \label{eq:upper-bound4}
\end{eqnarray}
where ${L}_1>0$ is a constant. 

Recall our assumptions that $0<\alpha <m/\ell$ where $\ell\in\{1,2\}$ and $\alpha =(1-q)m$, and also that $r_c(f)>(\ell m\alpha)/(m-\ell\alpha)$.
Setting $s=b$ in Lemma~\ref{lem:penrose}, we see that the left hand side of \eqref{eq:upper-bound4} is finite, so the first term on the right hand side of \eqref{eq:upper-bound2} is bounded by a constant which is independent of $N$.
For a non-negative random variable $Z>0$, we know that
\[
\mathbb{E}(Z) = \int_{0}^{\infty}\mathbb{P}(Z>z)\,dz,
\]
so the second term in \eqref{eq:upper-bound2} is bounded by
\begin{eqnarray}
\lefteqn{W_k\int_0^\infty\mathbb{P}\left(\max_{1\leq i\leq N}\|X_i\|^b > u\cdot N^{1-b/m}\right)\,du}\nonumber \\
	&& \quad \leq W_k\left[1 + N\int_1^\infty\mathbb{P}\left(\|X_1\|^b > \left(u^{m/(m-b)}N\right)^{1-b/m}\right)\,du\right] \nonumber\\
	\label{eq:upper-bound5}
\end{eqnarray}
%
By the Markov inequality $\mathbb{P}(Z>a)\leq\frac{1}{a}\mathbb{E}|Z|$ for $a>0$, we get for $u\geq 1$ that
\begin{eqnarray}
\lefteqn{\mathbb{P}\left(\|X_1\|^b > \left(u^{m/(m-b)}N\right)^{1-b/m}\right)} \qquad\qquad\qquad \nonumber\\
	& = & \mathbb{P}\left(\|X_1\|^{mb/(m-b)} > u^{m/(m-b)}N\right) \nonumber \\
	&\leq & \mathbb{E}\|X_1\|^{mb/(m-b)}\frac{1}{u^{m/(m-b)}N}.\label{eq:markov}
\end{eqnarray}
From \eqref{eq:upper-bound5} and \eqref{eq:markov},  we see that the second term in \eqref{eq:upper-bound2} is bounded by
\[
W_k\left[1+\int_1^\infty\mathbb{E}\|X_1\|^{mb/(m-b)}\frac{1}{u^{m/(m-b)}}\,du\right]
\]
which is independent of $N$, because we can choose $p$ to ensure that $0<1-b/m<1$, and 
\[
\mathbb{E}\|X_1\|^{\frac{mp(1-q)}{1-p(1-q)}} < \infty,
\text{ or equivalently }
r_{c}(f) > \frac{mp(1-q)}{1-p(1-q)},
\]
which is consistent with conditions of Theorem 2.1.
Note that the function $h(p,q)=\frac{mp(1-q)}{1-p(1-q)}$ is such that $h(1,q)$ gives the right-hand side of \eqref{eq:condition-conv-in-mean} and $h(2,q)$ gives the right-hand side of \eqref{eq:condition-conv-in-mean-square}. Moreover, if $r_c(f) > h(1,q)$ for some $q<1$ (resp. $r_c(f)>h(2,q)$ for some $q$ satisfying $1/2<q<1$), we also have $r_c(f)>h(p,q)$ for some $p>1$ (resp. $r_c(f)>h(p,q)$ for $p>2)$.

\section{Hypothesis tests}
\label{sec:htest}

We now restrict the class $\mathcal{K}$ to only those distributions which satisfy the following conditions: for any fixed $k\geq 1$ and $q>1/2$,
\begin{align*}
& \mathbb{E}(\hat{H}_{N,k,q}) \to H_q \quad\text{as $N\to\infty$, and} \\[0.5ex]
& \hat{H}_{N,k,q}\to H_q 				\quad\text{in probability as $N\to\infty$.}
\end{align*}
By Theorem~\ref{thm:main}, we know that $\mathcal{K}$ contains $T_{m}(a,\Sigma,\nu)$ for all $\nu>2$ and $P_m(a,\Sigma,\gamma)$ for all $\gamma>0$.

\bigskip
Let $X_1,X_2,\ldots,X_N$ be independent and identically distributed random vectors with common density $f\in\mathcal{K}$, and let $\hat{C}_N$ be the sample covariance matrix,
\[
\hat{C}_N = \frac{1}{N-1}\sum_{i=1}^N(X_i-\bar{X})(X_i-\bar{X})^{\prime}.
\]

\begin{enumerate}
\item
To test the hypothesis $X\sim T_m(a,\Sigma,\nu_0)$ where $\nu_0>2$, we define the test statistic
\begin{equation}\label{eq:tstat-student}
W^S_{N,k}(m,\nu) =  H^{\rm max}_q - \hat{H}_{N,k}(m,q),
\end{equation}
where
$H^{\rm max}_q = \frac{1}{2}\log|\hat{\Sigma}_N| + c'_S(m,q,\nu)$
with $q=1-1/(\nu+m)$ and $\displaystyle\hat{\Sigma}_N = (1-2/\nu)\hat{C}_N$.
\smallskip
\item
To test the hypothesis $X\sim P_m(a,\Sigma,\gamma_0)$ where $\gamma_0>0$, we define the test statistic
\begin{equation}\label{eq:tstat-pearsonII}
W^P_{N,k}(m,\gamma) =  H^{\rm max}_q - \hat{H}_{N,k}(m,q),
\end{equation}
where
$H^{\rm max}_q =  \frac{1}{2}\log|\hat{\Sigma}_N| + c'_P(m,q,\gamma)$
with $q = 1+1/\gamma$ and $\displaystyle\hat{\Sigma}_N = (2\gamma+m+2)\hat{C}_N$.
\end{enumerate}
By the law of a large numbers, $\hat{C}_N\to C$ in probability as $N\to\infty$, so by Slutsky's theorem, for any fixed $k\geq 1$, we have that
\[
\lim_{N\to\infty}W^S_{N,k}(m,\nu)\stackrel{P}\longrightarrow \begin{cases}
    0 & \text{if $X\sim T_m(a,\Sigma,\nu)$,} \\
    c > 0 & \text{otherwise},
\end{cases}
\]
and
\[
\lim_{N\to\infty}W^P_{N,k}(m,\gamma)\stackrel{P}\longrightarrow \begin{cases}
    0 & \text{if $X\sim P_m(a,\Sigma,\gamma)$,} \\
    c > 0 & \text{otherwise},
\end{cases}
\]
where ``$\stackrel{P}\longrightarrow$'' denotes convergence in probability and $c$ is a constant that depends on the distribution of $X$.

The distributions of $W^P_{N,k}(m,\nu)$ when $X\sim T_m(a,\Sigma,\nu)$ and $W^P_{N,k}(m,\gamma)$ when $X\sim P_m(a,\Sigma,\gamma)$ are unknown. An analytical derivation of these distributions seems difficult, because the random variables $\hat{H}_{N,k}$ and $\hat{C}_N$ are not independent and their covariance appears to be intractable, despite the fact that the asymptotic distribution of $\hat{H}_{N,k}$ can be revealed by applying the results of 
\textcite{chatterjee2008},
\textcite{penrose2011},
\textcite{delattre2017} and
\textcite{berrett2019}, 
and that of $\hat{C}_N$ by the delta method. In the next section, we investigate these null distributions using Monte Carlo methods.

\section{Numerical experiments}
\label{sec:numer}

\subsection{Random samples}
\label{subsec:random-samples}

Random samples from $T_m(a,\Sigma,\nu)$ and $P_m(a,\Sigma,\gamma)$ can be generated according to the stochastic representation
\[
X = RBU + a,
\]
where $R$ represents the radial distance $\big[(X-a)'\Sigma^{-1}(X-a)\big]^{1/2}$, $B$ is an $m\times m$ matrix with $B^TB=\Sigma$ and $U$ is uniformly distributed on the unit $m$-sphere $S^{m-1}$. In particular,
\begin{align*}
& R^2\sim\text{InvGamma}(m/2,m/2) \text{ yields } X\sim T_m(a,\Sigma,\nu), \text{ and}\\
& R^2\sim\text{Beta}(m/2,\gamma+1) \text{ yields } X\sim P_m(a,\Sigma,\gamma).
\end{align*}
%
%
Let $I_m$ be the $m\times m$ identity matrix.
We investigate the distributions 
\begin{align*}
& T_m(\nu) = T_m\big(0,I_m,\,\nu\big) \text{ for }\nu>2 \text{ and } \\
& P_m(\gamma) = P_m\big(0,I_m,\,\gamma\big) \text{ for }\gamma>1.
\end{align*}

\begin{figure*}[htb]
\begin{subfigure}[b]{.5\linewidth}
\centering\includegraphics[scale=0.45]{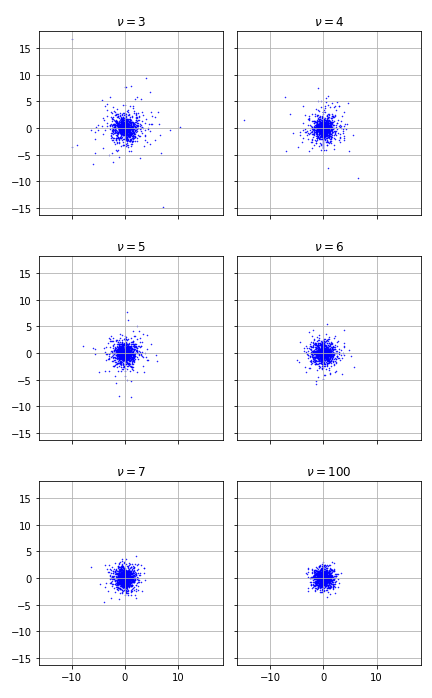}
\caption{Scatter plots for $T_2(\nu)$}\label{fig:scatter-MST-3x2}
\end{subfigure}
\hfill 
\begin{subfigure}[b]{.5\linewidth}
\centering\includegraphics[scale=0.45]{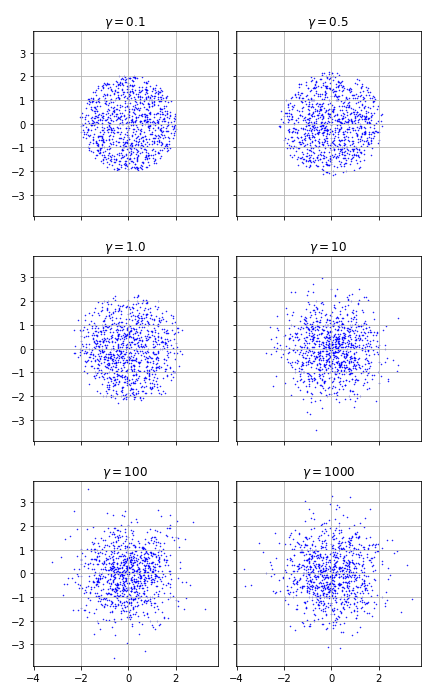}
\caption{Scatter plots for $P_2(\gamma)$}\label{fig:scatter-MP2-3x2}
\end{subfigure}%
\caption{Scatter plots for the bivariate Student and Pearson II distributions.}
\label{fig:scatter}
\end{figure*}

\subsection{Consistency}
\label{subsec:consistency}

\begin{figure*}[htb]
\centering\includegraphics[scale=0.45]{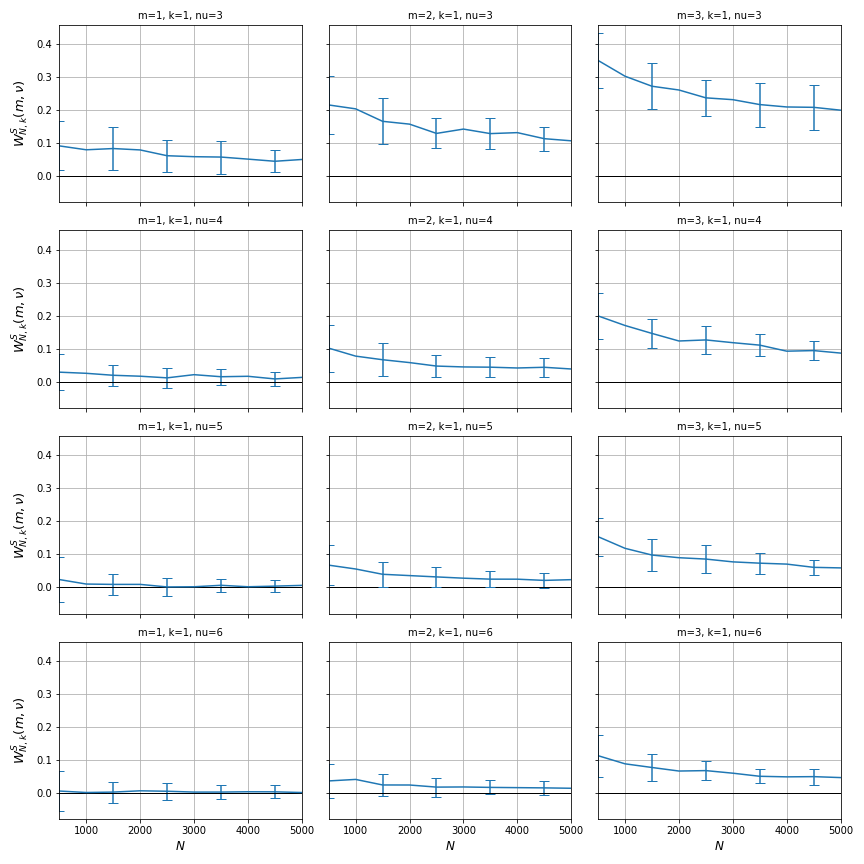}
\caption{The asymptotic behaviour of $W^S_{N,k}(m,\nu)$ as $N\to\infty$ for $k=1$. The rate of convergence appears to increase with $\nu$ and decrease with $m$.}\label{fig:cons-mst-k1}
\end{figure*}

\begin{figure*}[htb]
\centering\includegraphics[scale=0.55]{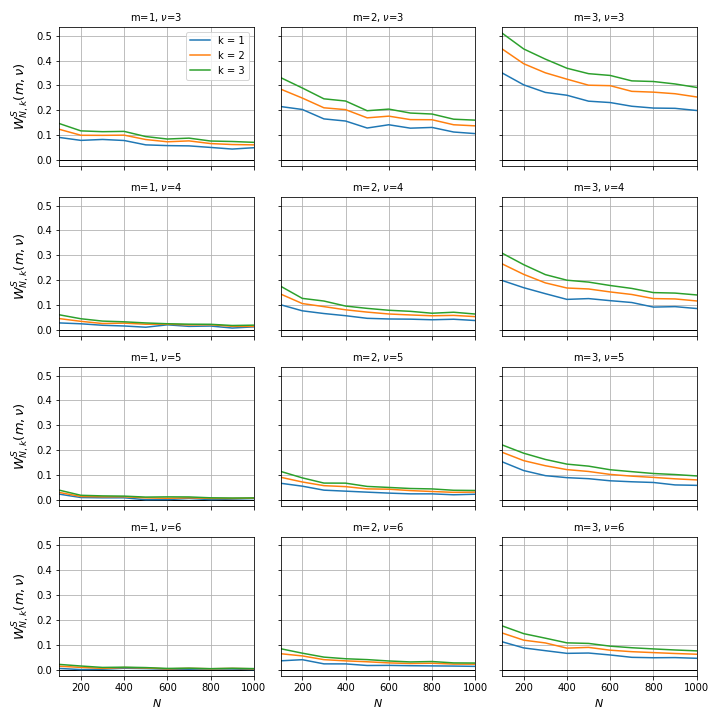}
\caption{The asymptotic behaviour of $W^S_{N,k}(m,\nu)$ as $N\to\infty$. The rate of convergence appears to increase with $\nu$ and decrease with $m$.}\label{fig:cons-mst-kall}
\end{figure*}

To investigate the consistency of $W^S_{N,k}(m,\nu)$ for various values of $m$ and $\nu$, we generate $M=100$ random samples of size $N$ from the $T_m(\nu)$ distribution, with $N$ increasing from $N=500$ to $N=5000$ in steps of $500$, and record the value of $W^S_{N,k}(m,\nu)$ for $k=1,2,3$ at each step. The mean values of the statistics for $k=1$ are shown in Figure~\ref{fig:cons-mst-k1}, where the lengths of the error bars are equal to the standard deviations of the statistics around their mean values. The mean statistics for $k=1,2,3$ are shown in Figure~\ref{fig:cons-mst-kall}, where it is evident that the rate of convergence increases with the parameter $\nu$ and decreases with the dimension $m$.

\begin{figure*}[htb]
\centering\includegraphics[scale=0.45]{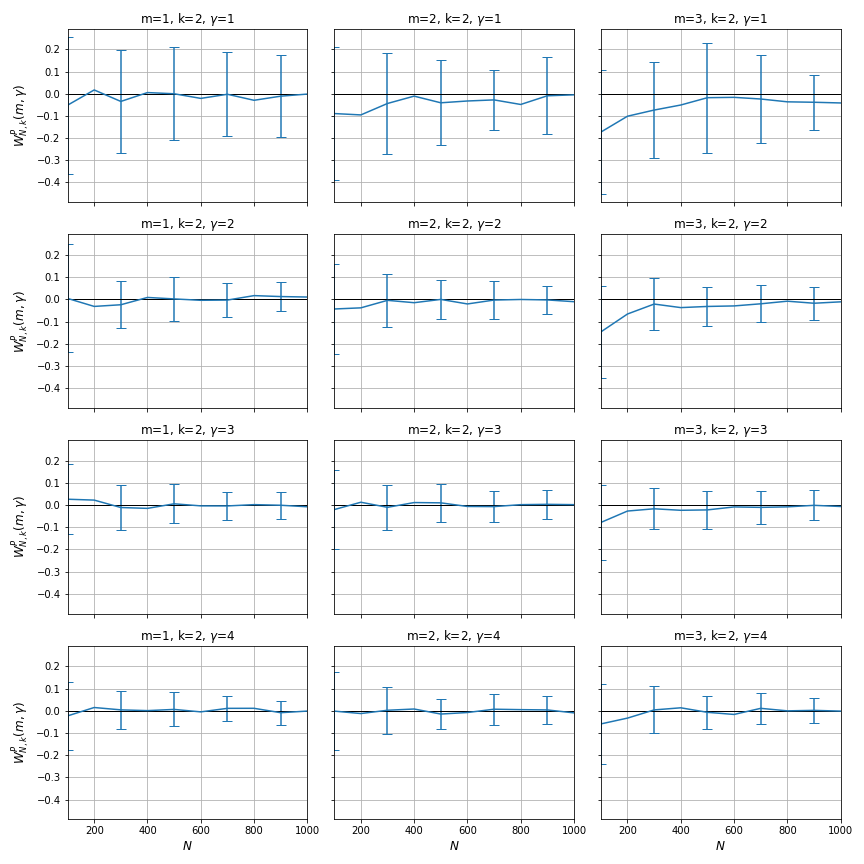}
\caption{Asymptotic behaviour of $W^P_{N,k}(m,\gamma)$ as $N\to\infty$ for $k=2$.}
\label{fig:cons-mp2-k2}
\end{figure*}

\begin{figure*}[htb]
\centering\includegraphics[scale=0.55]{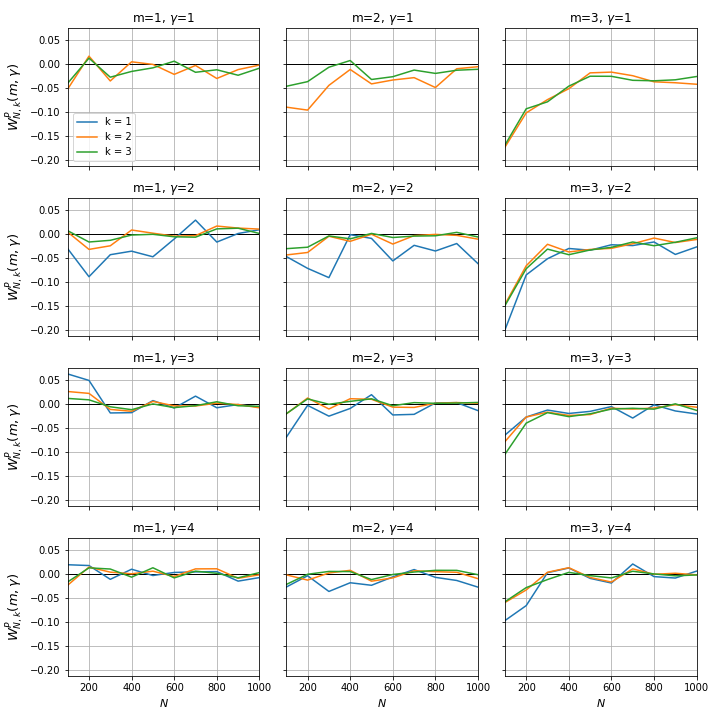}
\caption{Asymptotic behaviour of $W^P_{N,k}(m,\gamma)$ as $N\to\infty$. Note that the statistic is defined only for $k>1/\gamma$.}
\label{fig:cons-mp2-kall}
\end{figure*}

The experiment is repeated for $W^P_{N,k}(m,\gamma)$ but this time with samples increasing in size from $N=50$ to $N=500$ in steps of $50$. The mean values of the statistics for $k=2$ are shown in Figure~\ref{fig:cons-mp2-k2}, where lengths of the error bars are equal to the standard deviations of the statistics around their mean values. The mean statistics for $k=1,2,3$ are shown in Figure~\ref{fig:cons-mp2-kall}: note that these are only defined for $k>1/\gamma$. The convergence of $W^P_{N,k}(m,\gamma)$ is evidently much faster than that of $W^S_{N,k}(m,\nu)$, perhaps because the support of $P_m(\gamma)$ is bounded for any finite $\gamma>0$ while the support of $T_m(\nu)$ is unbounded.

\subsubsection{Rates of convergence}
\label{subsec:rates-of-convergence}

\begin{figure*}[htb]
\centering\includegraphics[scale=0.55]{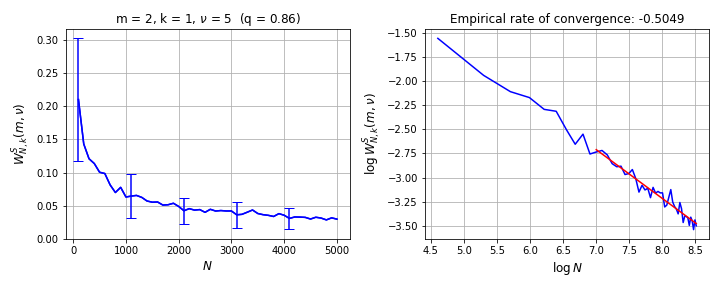}
\caption{Asymptotic behaviour of $W^S_{N,k}(m,\nu)$ with $m=2$, $k=1$ and $\nu=5$.}\label{fig:cons-nu5}
\end{figure*}

In Figure~\ref{fig:cons-nu5}, we plot the convergence of $W^S_{N,k}(m,\nu)$ as $N\to\infty$ with $m=2$, $k=1$ and $\nu=5$, together with the corresponding plot of $\log W^S_{N,k}(m,\nu)$ against $\log N$. THe latter suggests an empirical convergence rate of approximately $O(N^{-1/2})$ as $N\to\infty$. 

\begin{figure*}[htb]
\centering\includegraphics[scale=0.55]{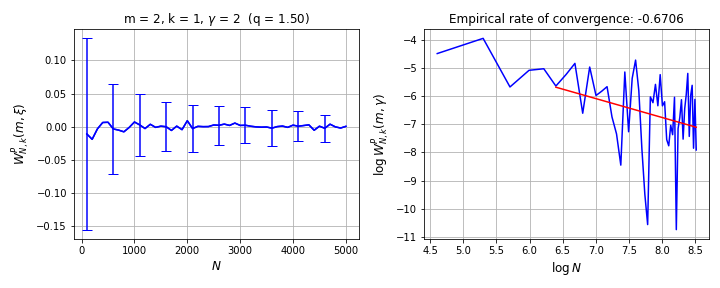}
\caption{Asymptotic behaviour of $W^P_{N,k}(m,\gamma)$ with $m=2$, $k=1$ and $\gamma=2$.}\label{fig:cons-xi2}
\end{figure*}

The experiment is repeated for $W^P_{N,k}(m,\gamma)$ with $m=2$, $k=2$ and $\gamma=2$. The results are shown in Figure~\ref{fig:cons-xi2}, which in this case suggest an empirical convergence rate of approximately $O(N^{-2/3})$ as $N\to\infty$. Analytic rates of convergence for $W^S_{N,k}(m,\nu)$ and $W^P_{N,k}(m,\gamma)$ are currently under investigation by the authors.

\subsection{Empirical distribution of the test statistics}
\label{subsec:empirical-dist}

\begin{figure*}[htb]
\centering\includegraphics[scale=.5]{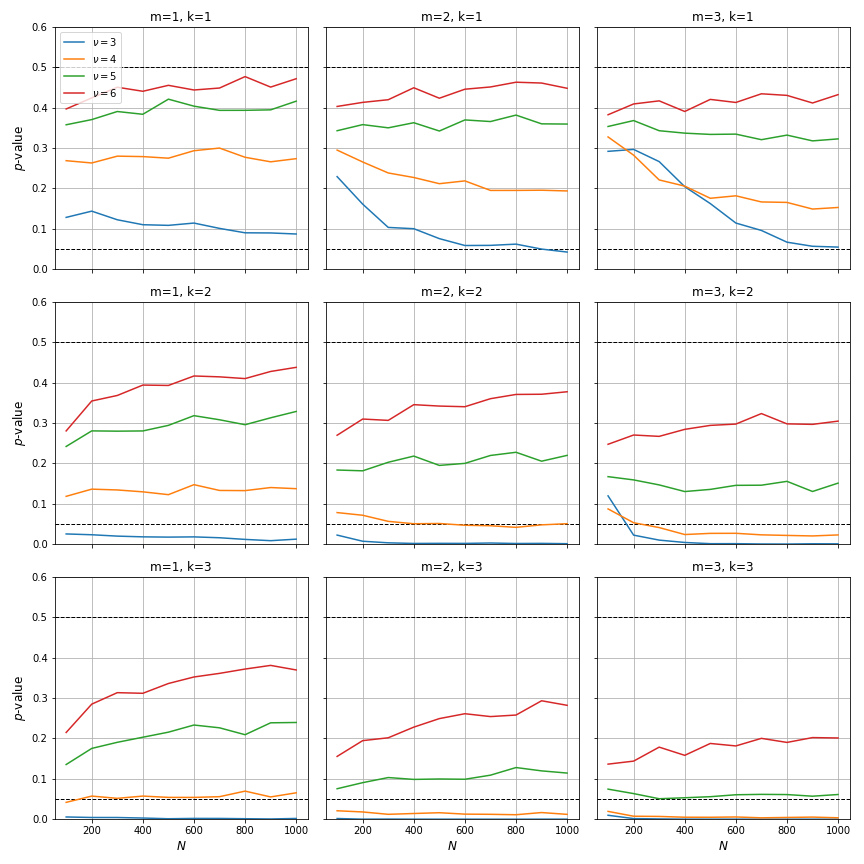}
\caption{Average Shapiro-Wilk probability values for $W^S_{N,k}(m,\nu)$ as $N$ increases for different values of $m$, $k$ and $\nu$ ($100$ repetitions).}
\label{fig:MST-sw-pvals-all}
\end{figure*}

For different values of $(N,k)$ and $(m,\nu)$, we generate $n=100$ random samples of size $N$ from the $T_m(\nu)$ distribution and record the value of $W^S_{N,k}(m,\nu)$ each time. We then apply the Shapiro-Wilk test for normality \parencite{shapiro1965} to this random sample and record the probability value computed by the test. This process is repeated $M=1000$ times. Figure~\ref{fig:MST-sw-pvals-all} illustrates how the mean probability value behaves as $N$ increases for various values of $m$, $k$ and $\nu$, where it appears that normal approximation improves as the distribution parameter $\nu$ increases, but deteriorates as $k$ increases.

\begin{figure*}[htb]
\centering\includegraphics[scale=.5]{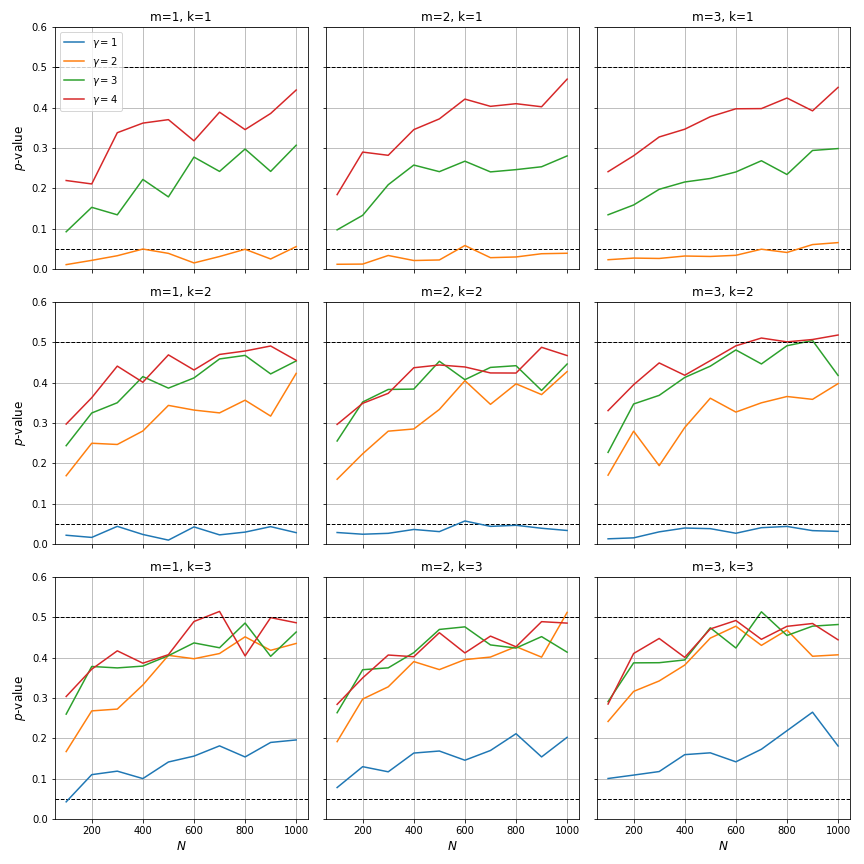}
\caption{Average Shapiro-Wilk probability values for $W^P_{N,k}(m,\gamma)$ as $N$ increases for different values of $m$, $k$ and $\gamma$ ($100$ repetitions). Note that the statistic is only defined for $k>1/\gamma$.}
\label{fig:MP2-sw-pvals-all}
\end{figure*}

The experiment is repeated for the null distribution of $W^P_{N,k}(m,\gamma)$ with the results shown in  Figure~\ref{fig:MP2-sw-pvals-all}, where it appears that normal approximation again improves as the distribution parameter $\gamma$ increases, but in this case appears to also improve as $k$ increases.

\subsection{Point estimation}
\label{subsec:pointest}

Point estimates for $\nu$ and $\gamma$ can be computed according to
\begin{align*}
\hat{\nu} & = {\rm argmin}_{\nu>2} W^S_{N,k}(m,\nu), \text{ and} \\
\hat{\gamma} & = {\rm argmin}_{\gamma>1} W^P_{N,k}(m,\gamma) \text{ respectively.}
\end{align*}

\begin{figure*}[htb]
\centering\includegraphics[scale=0.45]{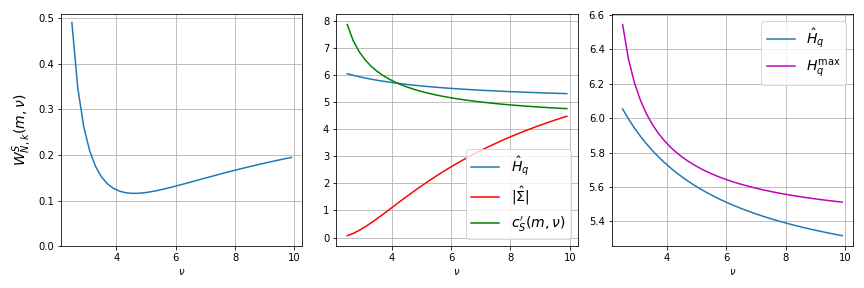}
\caption{$W^S_{N,k}(m,\nu)$ with $N=1000$, $k=1$, $m=3$ and $\nu_0=4$ ($q_0=0.86$).\\ 
The point estimate is $\hat{\nu}=4.7$.}\label{fig:pest-nu4}
\end{figure*}

A random sample of size $N=1000$ was generated from the $T_3(\nu)$ distribution with $\nu=4$, and the value of $W^S_{1000,1}(3,\nu)$ was then computed for different values of $\nu$ in the range $[2.5,10]$. The results are shown in Figure~\ref{fig:pest-nu4}, where we see that the statistic reaches a minimum value at approximately $\nu=4$. Note that because we take $\Sigma=I_m$, the estimated determinant $|\hat{\Sigma}|$ is approximately equal to $1$ when $\nu=4$. 

\begin{figure*}[htb]
\centering\includegraphics[scale=0.45]{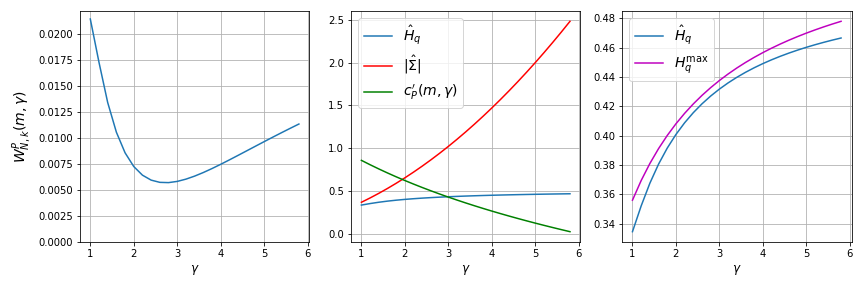}
\caption{$W^P_{N,k}(m,\gamma)$ with $N=1000$, $k=1$, $m=3$ and $\gamma_0=3$ ($q_0=1.33$).\\
The point estimate is $\hat{\gamma}=2.8$.}\label{fig:pest-gamma3}
\end{figure*}

The experiment is repeated for a random sample from the $P_3(\gamma)$ distribution with $\gamma=3$, and the value of $W^P_{1000,1}(3,\gamma)$ computed for different values of $\gamma$ in the range $[1,6]$. The results are shown in Figure~\ref{fig:pest-gamma3}, where we see that the statistic reaches a minimum value at approximately $\gamma=3$. 

The theoretical properties of these estimators are currently under investigation by the authors.

\bibliographystyle{plain}  
\bibliography{entropy} 


\end{document}